\documentclass[a4paper]{article}
\usepackage{amsmath,latexsym,epsfig,graphics,amsfonts,amssymb,mathrsfs,amsthm}

\newtheorem{theorem}{Theorem}[section]
\newtheorem{lemma}[theorem]{Lemma}
\newtheorem{proposition}[theorem]{Proposition}
\newtheorem{corollary}[theorem]{Corollary}
\theoremstyle{definition}

\theoremstyle{plain}

\numberwithin{equation}{section}
\numberwithin{figure}{section}

\renewcommand{\leq}{\leqslant}
\renewcommand{\geq}{\geqslant}

\begin{document}

\title{Ford circles, continued fractions, and best approximation of the second kind}
\author{Ian Short \\
Centre for Mathematical Sciences\\ 
Wilberforce Road\\
Cambridge CB3 0WB \\
United Kingdom\\
\texttt{ims25@cam.ac.uk}}
\date{\today}

\maketitle 



\begin{abstract}
We give an elementary geometric proof using Ford circles that the convergents of the continued fraction expansion of a real number $\alpha$ coincide with the rationals that are best approximations of the second kind of $\alpha$.
\end{abstract}

\section{Introduction}\label{I}

This paper is about a geometric view of the relationship between continued fractions and approximation of real numbers by rationals. Whenever we speak of a rational $u/v$ we mean that $u$ and $v$ are coprime integers and $v$ is positive.  Given a real number $\alpha$ we follow Khinchin \cite[Section 6]{Kh1997} in describing a rational $a/b$ as a \emph{best approximation of the second kind} of $\alpha$ provided that, for each rational $c/d$ such that $d\leq b$, we have
\[
\left| b\alpha-a\right| \leq \left| d\alpha-c\right|,
\]
with equality if and only if $c/d = a/b$.
Khinchin also defines \emph{best approximation of the first kind}, but that concept does not concern us here.

A \emph{continued fraction} is an expression of the form
\begin{equation*}
               b_0+     \cfrac{1}{b_1+
                      \cfrac{1}{b_2+
                       \cfrac{1}{b_3+ \dotsb
}}}\;,
\end{equation*}
where $b_0$ is an integer and the other $b_i$ are positive integers. Either the sequence $b_0,b_1,b_2,\dotsc$ is infinite, in which case the continued fraction is said to be \emph{infinite}, or there is a final member $b_N$ of this sequence, in which case the continued fraction is said to be \emph{finite}. Each of our finite continued fractions with $N\geq 1$ is assumed to satisfy $b_N\geq 2$. The \emph{convergents} of $\alpha$ are the rationals $A_n/B_n$ where,
\[
\frac{A_0}{B_0} = b_0\,,\qquad \frac{A_1}{B_1} = b_0+\frac{1}{b_1},\qquad \frac{A_2}{B_2}=b_0+\cfrac{1}{b_1+\cfrac{1}{b_2+}}\,,\,\dotsb.
\]
The \emph{value} of a finite continued fraction is the value of the final convergent $A_N/B_N$ (which is a rational number) and the \emph{value} of an infinite continued fraction is the limit of the sequence $A_n/B_n$ (which is an irrational number). To each real number $\alpha$ there corresponds a unique continued fraction with value $\alpha$.

All these facts about continued fractions can be found in \cite{Kh1997}, as can the next theorem (\cite[Theorems 16 and 17]{Kh1997}).

\begin{theorem}\label{H}
A rational $x$, which is not an integer, is a convergent of a real number $\alpha$ if and only if it is a best approximation of the second kind of $\alpha$.
\end{theorem}

It is convenient to assume that $x$ is not an integer in Theorem~\ref{H}, and later on, to avoid tiresome discussions of this trivial case. Theorem~\ref{H} fails when $x$ is an integer, because if $m+\tfrac12 \leq \alpha < m+1$, for some integer $m$, then $m$ is a convergent of $\alpha$, but not a best approximation of the second kind of $\alpha$. A version of Theorem~\ref{H} including the possibility that $x$ is an integer can be found in \cite[Theorem 1]{Ir1989}.


Classic proofs of Theorem~\ref{H}, such as that given in \cite{Kh1997}, are algebraic. Irwin proves Theorem~\ref{H} using plane lattices in \cite{Ir1989}. Our aim is to give an illuminating geometric proof based on the theory of Ford circles. Ford circles, developed by Ford in \cite{Fo1938}, are objects most naturally associated with hyperbolic geometry, and our proof has undertones of hyperbolic geometry. We now give a brief description of Ford circles and their relationship to continued fractions (full details can be found in \cite{Fo1938}).

\begin{figure}[ht]
\centering
\includegraphics[scale=1.0]{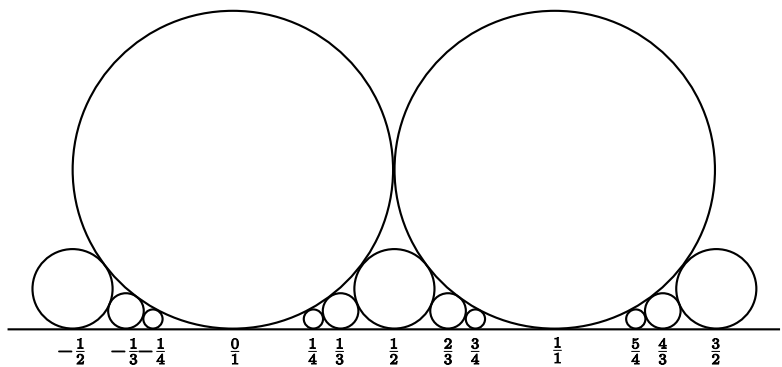}
\caption{Ford circles.}
\label{S}
\end{figure}

The \emph{Ford circle} $C_x$ of a rational number $x=a/b$ is the circle in the complex plane with centre $x+i/(2b^2)$ and radius $1/(2b^2)$. This circle is tangent to the real axis at $x$, and otherwise lies in the upper half-plane. A selection of Ford circles are shown in Figure~\ref{S}.

\begin{figure}[ht]
\centering
 \includegraphics[scale=1.0]{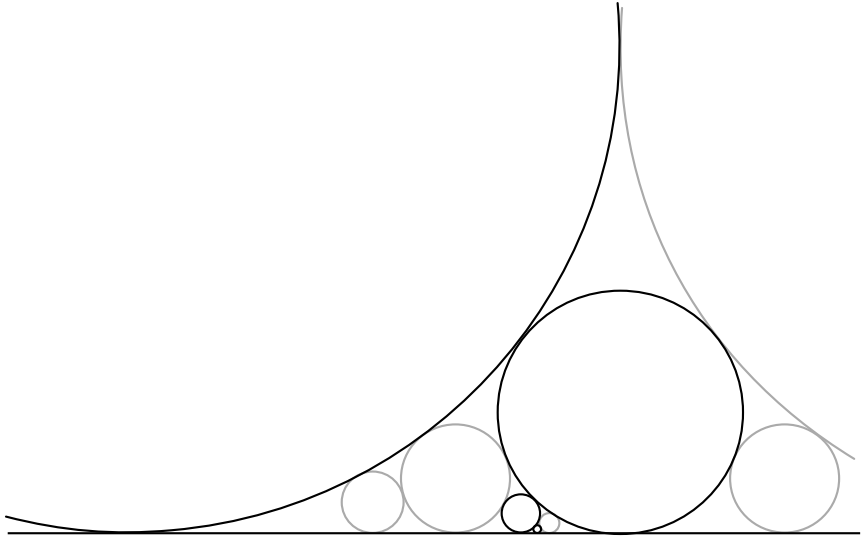}
\caption{A continued fraction chain.}
\label{T}
\end{figure}

Two Ford circles $C_x$ and $C_y$, where $x=a/b$ and $y=c/d$, are tangent if and only if $|ad-bc|=1$, and if they are not tangent then they are wholly external to one another (Ford circles do not overlap). We define the \emph{continued fraction chain} of a real number $\alpha$ to be the sequence of Ford circles $C_{A_0/B_0}, C_{A_1/B_1}, C_{A_2/B_2},\dotsc$, where $A_n/B_n$ are the convergents of $\alpha$. Since $|A_nB_{n-1}-A_{n-1}B_n|=1$ we see that each pair of consecutive circles in the continued fraction chain of $\alpha$ are tangent. Also, the $B_n$ are given by the recurrence relation $B_0=1$, $B_1=b_1$, and $B_n=b_nB_{n-1}+B_{n-2}$ for $n\geq 2$, which means that the sequence $B_1,B_2,\dotsc$ of positive integers is increasing. Consequently, the sequence of radii $1/(2B_1^2), 1/(2B_2^2),\dotsc$ of the circles from the continued fraction chain is decreasing. Finally, the members of the continued fraction chain alternate from the left to the right side of $\alpha$, because
\begin{equation}\label{W}
\frac{A_0}{B_0} < \frac{A_2}{B_2} < \frac{A_4}{B_4} <\dotsb < \alpha < \dotsb < \frac{A_5}{B_5} < \frac{A_3}{B_3} < \frac{A_1}{B_1}.
\end{equation}
The first few Ford circles from a continued fraction chain are shown in Figure~\ref{T} (in black).

A circle $C$ that is tangent to the real axis at a point $z$, and that otherwise lies in the upper half-plane, is called a \emph{horocircle}. We denote the radius of $C$ by  $\textnormal{rad}[C]$, and describe the point $z$ as the \emph{base point} of $C$.  In order to state our geometric version of Theorem~\ref{H}, we introduce a new definition. In this definition we use the simple fact that, given a real number $z$ and a horocircle $C$, there is a unique horocircle $D$ that is tangent to $C$ and has base point $z$. If $C$ has base point $z$ then we consider $D$ to have radius $0$. 

Given a real number $\alpha$ we say that a Ford circle $C_x$ is \emph{nearby} to $\alpha$ if, for each Ford circle $C_z$ other than $C_x$ with $\textnormal{rad}[C_z]\geq \textnormal{rad}[C_x]$, the radius of the unique horocircle tangent to $C_z$ and with base point $\alpha$ is larger than the radius of the unique horocircle tangent to $C_x$ and with base point $\alpha$. When $\alpha$ is rational, $C_\alpha$ is nearby to $\alpha$, but no Ford circle with equal or smaller radius than $C_\alpha$ is nearby to $\alpha$.

\begin{theorem}\label{U}
Let $\alpha$ be a real number. Given a rational $x$, which is not an integer, the following are equivalent:
\begin{enumerate}
\item[\textnormal{(i)}] $x$ is a convergent of $\alpha$;
\item[\textnormal{(ii)}] $C_x$ is a member of the continued fraction chain of $\alpha$;
\item[\textnormal{(iii)}] $x$ is a best approximation of the second kind of $\alpha$;
\item[\textnormal{(iv)}] $C_x$ is nearby to $\alpha$;
\item[\textnormal{(v)}] there is a Ford circle $C_y$ tangent to $C_x$ such that $\textnormal{rad}[C_x] > \textnormal{rad}[C_y]$, and  either $\alpha=x$ or $\alpha$ lies in the open interval bounded by $x$ and $y$.
\end{enumerate}
\end{theorem}

\begin{figure}[ht]
\centering
\includegraphics[scale=1.0]{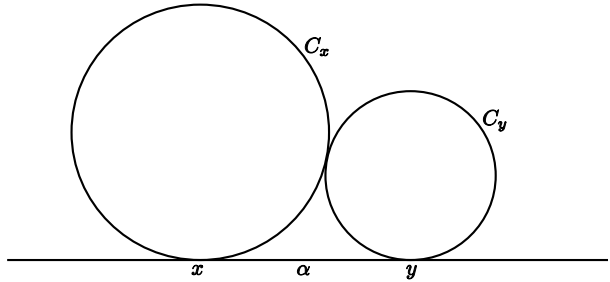}
\caption{The geometry of Theorem~\ref{U} (v).}
\label{Z}
\end{figure}

Statement (v) is illustrated in Figure~\ref{Z}.

Statement (ii) of Theorem~\ref{U} is merely a geometric reformulation of statement (i), and in the next section we see that statement (iv) is a geometric reformulation of statement (iii). The equivalence of statements (i) and (iii) yields Theorem~\ref{H}.

\section{Best approximation of the second kind}\label{K}

The key idea in this paper is about explaining best approximation of the second kind in terms of Ford circles so that, using Ford's continued fraction chains, we can prove Theorem~\ref{H} geometrically. Our key idea is encapsulated in the next proposition.

\begin{proposition}\label{D}
Given a real number $\alpha$, the rational $x$ is a best approximation of the second kind of $\alpha$ if and only if $C_{x}$ is nearby to $\alpha$.
\end{proposition}

To prove Proposition~\ref{D} we need the next lemma and corollary.

\begin{lemma}\label{A}
Two horocircles $C$ and $D$ with radii $r$ and $s$ and distinct base points $x$ and $y$, which intersect in at most one point, satisfy $|x-y|^2\geq 4rs$, with equality if and only if $C$ and $D$ are tangent.
\end{lemma}
\begin{figure}[ht]
\centering
\includegraphics[scale=1.0]{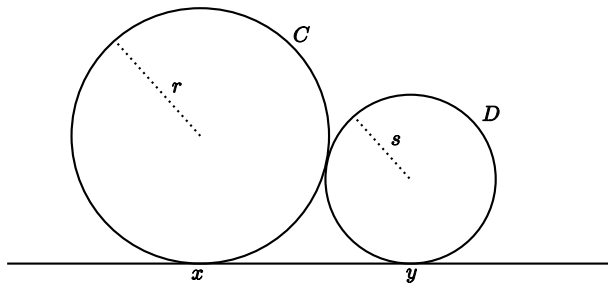}
\caption{The geometry of Lemma~\ref{A}.}
\label{C}
\end{figure}
\begin{proof}
Let $d$ be the distance between the centres of the two horocircles. Then $d\geq r+s$, with equality if and only if $C$ and $D$ are tangent. We can calculate $d$ by applying Pythagoras's Theorem to the triangle with vertices $x+ir$, $y+is$, and $x+is$, and the result follows immediately.
\end{proof}

\begin{corollary}\label{B}
The radius of the horocircle that is tangent to the Ford circle $C_{a/b}$, and has base point $\alpha$, is 
\[
\tfrac12 |b\alpha -a|^2.
\] 
\end{corollary}
\begin{proof}
This corollary follows from Lemma~\ref{A}, because $C_{a/b}$ has radius $1/(2b^2)$.
\end{proof}

\begin{proof}[Proof of Proposition~\ref{D}]
For each rational $z$, let $D_z$ denote the unique horocircle with base point $\alpha$ that is tangent to $C_z$. Now, a rational $x=a/b$ is a best approximation of the second kind of $\alpha$ if and only if for each rational $y=c/d$ distinct from $x$ and such that $d\leq b$ we have $|d\alpha-c|>|b\alpha-a|$. Equivalently, using Corollary~\ref{B}, for each Ford circle $C_y$ distinct from $C_x$ and such that $\textnormal{rad}[C_y]\geq \textnormal{rad}[C_x]$ we have $\textnormal{rad}[D_y]> \textnormal{rad}[D_x]$. In other words $x$ is a best approximation of the second kind of $\alpha$ if and only if $C_x$ is nearby to $\alpha$.
\end{proof}

\section{Ford circles}\label{R}

This section contains two elementary lemmas about basic properties of Ford circles. 

\begin{lemma}\label{X}
Let $C_x$ and $C_y$ be tangential Ford circles. If a rational $z$ lies strictly between $x$ and $y$ then $C_z$ has smaller radius than both $C_x$ and $C_y$.
\end{lemma}

\begin{figure}[ht]
\centering
\includegraphics[scale=1.0]{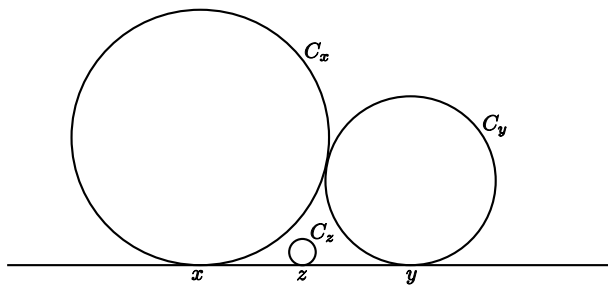}
\caption{The geometry of Lemma~\ref{X}.}
\label{N}
\end{figure}

\begin{proof}
Let $C_x$, $C_y$, and $C_z$ have radii $r_x$, $r_y$, and $r_z$. By Lemma~\ref{A}
\[
|x-y|^2=4r_xr_y,\quad |y-z|^2\geq 4r_yr_z, \quad |z-x|^2\geq 4r_zr_x.
\]
Hence
\[
r_z \leq \frac{|z-x|^2}{4r_x} < \frac{|x-y|^2}{4r_x} = r_y, 
\]
and similarly $r_z<r_x$.
\end{proof}

\begin{lemma}\label{Q}
Let $C_x$ and $C_y$ be tangential Ford circles such that $\textnormal{rad}[C_x]>\textnormal{rad}[C_y]$, and suppose that a real number $\alpha$ lies strictly between $x$ and $y$, and a rational $z$ lies strictly outside the interval bounded by $x$ and $y$. Then the radius of the horocircle that is tangent to $C_x$ and has base point $\alpha$ is smaller than the radius of the horocircle that is tangent to $C_z$ and has base point $\alpha$. 
\end{lemma} 

\begin{figure}[ht]
\centering
\includegraphics[scale=1.0]{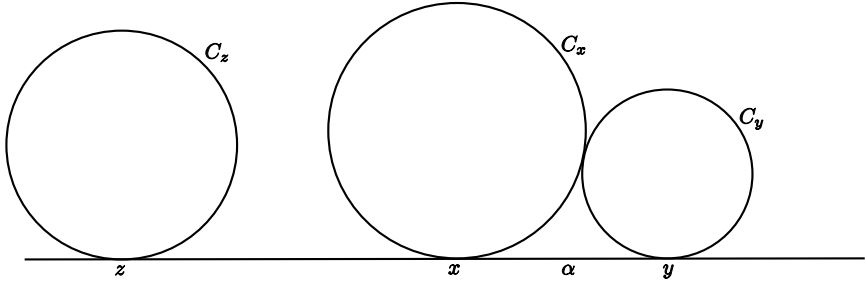}
\caption{The geometry of Lemma~\ref{Q}.}
\label{M}
\end{figure}

\begin{proof}
Let $C_x$, $C_y$, and $C_z$ have radii $r_x$, $r_y$, and $r_z$. Denote the radius of the horocircle that is tangent to $C_x$, and has base point $\alpha$, by $s_x$, and denote the radius of the horocircle that is tangent to $C_z$, and has base point $\alpha$, by $s_z$. By Lemma~\ref{A} we have
\[
|x-y|^2=4r_xr_y,\quad |x-\alpha|^2=4r_xs_x,\quad |x-\alpha|\leq |x-y|,
\]
from which it follows that
\begin{equation}\label{F}
s_x = \frac{|x-\alpha|^2}{4r_x} < \frac{|x-y|^2}{4r_x} = r_y.  
\end{equation}
By Lemma~\ref{A} we also have 
\[
|z-\alpha|^2=4r_zs_z,\quad |z-x|^2\geq 4r_xr_z,\quad |z-y|^2\geq 4r_yr_z.
\]
Depending on the order of $x$, $y$, and $z$ we either have 
\[
|z-\alpha|^2>|z-x|^2\geq 4r_xr_z > 4r_yr_z
\]
or 
\[
|z-\alpha|^2>|z-y|^2\geq 4r_yr_z.
\]
In both cases we obtain
\begin{equation}\label{G}
s_z = \frac{|z-\alpha|^2}{4r_z} >r_y.
\end{equation}
Combining \eqref{F} and \eqref{G} we conclude that $s_x < r_y<s_z$.
\end{proof} 

\begin{corollary}\label{Y}
Let $C_x$ and $C_y$ be tangential Ford circles such that $\textnormal{rad}[C_x]>\textnormal{rad}[C_y]$, and suppose that a real number $\alpha$ lies strictly between $x$ and $y$. Then $C_x$ is nearby to $\alpha$. 
\end{corollary}
\begin{proof}
Choose a Ford circle $C_z$ distinct from $C_x$ with $\textnormal{rad}[C_z]\geq \textnormal{rad}[C_x]$. By Lemma~\ref{X}, $z$ lies strictly outside the interval bounded by $x$ and $y$. By Lemma~\ref{Q}, the radius of the horosphere based at $\alpha$ and tangent to $C_x$ is smaller than the radius of the horosphere based at $\alpha$ and tangent to $C_z$. Hence $C_x$ is nearby to $\alpha$.
\end{proof}

\section{Proof of Theorem~\ref{U}}\label{J}
 
We can now prove Theorem~\ref{U}. In our proof we denote the convergents of $\alpha$ by $A_0/B_0, A_1/B_1, A_2/B_2,\dotsc$, and we define $r_n$ to be the radius $1/(2B_n^2)$ of $C_{A_n/B_n}$, so that $r_1, r_2,\dotsc$ is a strictly decreasing sequence. Statements (i) and (ii) of Theorem~\ref{U} are equivalent by the  definition of a continued fraction chain. Statements (iii) and (iv) are equivalent because of Proposition~\ref{D}. We proceed to prove that (i) implies (v), (v) implies (iv), and (iv) implies (i). 

It is convenient to first dismiss the two cases when $\alpha$ is rational, and $x$ is either the last or penultimate convergent of $\alpha$ (namely $A_N/B_N$ or $A_{N-1}/B_{N-1}$). If $x=A_N/B_N$ then (i) and (iv) are satisfied by definition, and (v) is also satisfied by choosing any Ford circle $C_y$  tangent to $C_x$ such that $\textnormal{rad}[C_y]<\textnormal{rad}[C_x]$. Suppose that $x=A_{N-1}/B_{N-1}$. Again, (i) holds by definition. Let $u=A_{N}-A_{N-1}$ and $v=B_{N}-B_{N-1}$. This pair are coprime because $|B_N(A_{N}-A_{N-1})-A_N(B_{N}-B_{N-1})|=1$, and because $B_N=b_NB_{N-1}+B_{N-2}$ and $b_N\geq 2$, we see that $B_{N-1}<v<B_N$. Let $y=u/v$. Notice that $\alpha=A_N/B_N$ lies strictly between $x$ and $y$. Therefore (v) is satisfied, and (iv) is satisfied by Corollary~\ref{Y}. Henceforth we assume that, when $\alpha$ is rational, $x\neq A_N/B_N$ and $x\neq A_{N-1}/B_{N-1}$. 

Now we show that (i) implies (v). Suppose that $x=A_n/B_n$ for some $n$. Define $y=A_{n+1}/B_{n+1}$. By \eqref{W}, $\alpha$ lies strictly between $x$ and $y$, and since $r_1,r_2,\dotsc$ is decreasing we have that $\textnormal{rad}[C_x]>\textnormal{rad}[C_y]$.

Next we show that (v) implies (iv). Suppose that (v) holds. The case $x=\alpha$ has already been dealt with, hence $\alpha$ lies strictly between $x$ and $y$, and it follows from Corollary~\ref{Y} that $C_x$ is nearby to $\alpha$.

Last we show that (iv) implies (i). Suppose that $C_x$ is nearby to $\alpha$. The decreasing sequence $r_1, r_2,\dotsc$ has limit $0$ if $\alpha$ is irrational, and limit $\textnormal{rad}[C_\alpha]$ if $\alpha$ is rational. In the latter case, as  $x=\alpha$ has been considered already, we have $\textnormal{rad}[C_x]>\textnormal{rad}[C_\alpha]$. Since $r_0=\tfrac12$ -- the greatest possible radius of a Ford circle -- either (a) $r_1=r_0$ and there is a unique integer $n\geq 1$ such that $r_n\geq \textnormal{rad}[C_x] > r_{n+1}$, or (b) $r_1<r_0$ and there is a unique integer $n\geq 0$ such that $r_n\geq \textnormal{rad}[C_x] > r_{n+1}$. Also, since $x$ is neither the last nor the penultimate convergent of $\alpha$, we see that $\alpha$ lies strictly between $A_n/B_n$ and $A_{n+1}/B_{n+1}$. If $x\neq A_n/B_n$ then, by Lemma~\ref{X}, $x$ lies outside the closed interval bounded by $A_n/B_n$ and $A_{n+1}/B_{n+1}$; however, this cannot be, because the assumption that $C_x$ is nearby to $\alpha$ then contradicts Lemma~\ref{Q}. Hence $x=A_n/B_n$ (and case (b) with $n=0$ cannot arise because $x$ is not an integer).

\section{Concluding remarks}\label{L}

Let $j$ denote the point $(0,0,1)$ in three-dimensional Euclidean space. The \emph{Ford sphere} $S_x$ of $x=a/b$, where $a$ and $b$ are coprime Gaussian integers, is the sphere with centre $x+j/(2|b|^2)$ and radius $1/(2|b|^2)$. Ford spheres share many properties with Ford circles, and they can be used in the study of Gaussian integer continued fraction expansions of complex numbers. A brief account can be found at the end of \cite{Fo1938}. It would be of interest to investigate whether the techniques of this paper can be applied to Ford spheres and Gaussian integer continued fractions to give results on approximation of complex numbers by quotients of Gaussian integers.


\end{document}